\documentclass[11pt]{article}

\usepackage[T1]{fontenc}
\usepackage[utf8]{inputenc}

\usepackage{amsmath,amssymb,amsthm,mathtools}
\usepackage{mathrsfs}
\usepackage{enumitem}
\usepackage{hyperref}
\usepackage{graphicx}
\usepackage{longtable}
\usepackage{xcolor}
\usepackage{listings}
\usepackage[all]{xy}

\title{From Matrices to Morphisms II: Examples of Computational Categories in MATLAB and Octave}
\author{
Nelson Martins-Ferreira\footnote{martins.ferreira@ulo.pt}\\
Universidade de Leiria e Oeste\\
(Polytechnic of Leiria)\\
Portugal
}
\date{\today}

\newtheorem{theorem}{Theorem}[section]

\theoremstyle{definition}

\newtheorem{example}[theorem]{Example}

\begin{document}

\maketitle

\begin{abstract}
This article develops the notion of a computational category, motivated by the observation that many mathematical structures occurring in scientific computing already admit concrete algorithmic representations. Finite sets, finite-dimensional vector spaces, relations, and graphs may all be represented through elementary MATLAB and Octave data structures such as indexing vectors, matrices, logical arrays, and sparse matrices.

The central idea is that a category may often be described through a category of canonical representatives indexed by the natural numbers. This viewpoint leads naturally to representative categories whose objects are dimensions and whose morphisms are finite computational data structures. The resulting framework unifies several examples, including matrix categories, index categories, categories of column spaces, row-based models of finite sets, and categories combining finite subsets with finite-dimensional vector spaces through mixed morphisms.

The paper also introduces computational categories with distinguished non-finite representatives, extending the finite setting while retaining concrete coordinate descriptions. Throughout, categorical constructions are interpreted algorithmically and implemented using standard MATLAB and Octave operations.

The work contributes to a computational perspective on category theory and, conversely, to a categorical understanding of matrix-oriented scientific computing.
\end{abstract}

 \medskip

\section{Introduction} 

Category theory provides a powerful language for describing mathematical structures and the relationships between them. Over the past decades, categorical ideas have influenced a wide range of areas including computer science, programming language design, formal methods, databases, and mathematical software \cite{Adamek,Awodey,Lawvere,Leinster,MacLane}.
 At the same time, environments such as MATLAB and Octave have become standard tools for scientific computing, numerical analysis, engineering, data processing, and computational experimentation \cite{Octave,Matlab,Moler}. Although these two traditions evolved largely independently, both are fundamentally concerned with composition, structure, and the representation of mathematical processes.

 The starting point of this paper is the observation that many computational structures routinely manipulated in MATLAB and Octave already possess natural categorical descriptions. Finite sets may be represented through indexing vectors, finite-dimensional vector spaces through matrices, relations through logical matrices, and graphs through sparse matrices. These representations are not merely convenient data structures: they carry morphisms, compositions, factorizations, limits, colimits, and other categorical constructions that can be implemented directly through standard computational operations.
 
  This observation motivates the notion of a \emph{computational category}. Informally, a computational category is a category whose objects and morphisms admit concrete computational representations through a distinguished collection of canonical representatives. In such a setting, abstract mathematical objects may be described through finite data structures, while categorical constructions become executable algorithms. The resulting viewpoint provides a bridge between categorical reasoning and scientific computing. 
  
  The present work develops this idea through a collection of examples arising naturally in matrix-oriented computation. We study categories of matrices and indexing vectors, categories of row and column representations, categories of relations and spans, and categories combining finite combinatorial structures with finite-dimensional linear algebra. A common theme throughout is the passage from abstract mathematical objects to representative objects indexed by dimensions, allowing morphisms to be encoded by matrices, indexing vectors, logical arrays, sparse matrices, and related computational data structures. 
  
  The paper should be viewed as a companion to the manuscript \emph{From Matrices to Morphisms I: Towards Categorical Programming with MATLAB and Octave} \cite{MFGTLab}, which develops a broader categorical interpretation of standard MATLAB and Octave constructions and builds upon earlier work in \cite{MFScripta2025}. Whereas that work focuses primarily on categorical constructions arising from existing commands, the present article focuses on the computational categories themselves, their representative forms, and the relationships between the various computational models that arise from finite sets, vector spaces, relations, graphs, and geometric data. 
  
  The philosophy adopted throughout is deliberately pragmatic. Rather than proposing a new programming language or a specialized category-theoretic software framework, we seek to identify categorical structures already present in ordinary computational practice. The examples developed in the paper suggest that many familiar MATLAB and Octave operations may be interpreted as concrete realizations of categorical constructions, and that several common computational data structures can themselves be understood as representative categories. 
  
  The paper is largely self-contained. For the convenience of readers coming from either mathematics or computing, the appendices include a brief introduction to Octave programming, a concise review of the categorical notions used in the text, supplementary computational material, and reference implementations of several constructions discussed throughout the paper.

\section{A Brief Introduction to Octave Programming}

The present notes assume a basic familiarity with Octave (or MATLAB). Nevertheless, a few remarks concerning the numerical model used by the system are useful.

Octave stores integers and real numbers using finite machine representations. In particular, ordinary numerical computations are performed using the IEEE 754 double-precision floating-point standard. As a consequence, not every mathematical identity remains valid when interpreted computationally.

For example, although integers up to approximately \(2^{53}\) may be represented exactly, one encounters phenomena such as

\begin{lstlisting}
a = 2^53;
isequal(a,(a+1)-1)
\end{lstlisting}
which returns \texttt{false}. On the other hand,

\begin{lstlisting}
a = 2^53-1;
isequal(a,(a+1)-1)
\end{lstlisting}
returns \texttt{true}. When exact integer arithmetic is required, one may instead use integer types such as

\begin{lstlisting}
a = int64(2^53);
isequal(a,(a+1)-1)
\end{lstlisting}

Throughout these notes, such technical issues will play no essential role. We shall simply regard the natural numbers as represented by machine integers within the available numerical range.

The same remarks apply to real numbers. Since computations are performed using double precision, equality should always be interpreted computationally rather than mathematically. For example,

\begin{lstlisting}
isequal(1,1+eps)
isequal(1,1+eps/2)
\end{lstlisting}
returns

\begin{verbatim}
0
1
\end{verbatim}
respectively, where

\[
\texttt{eps}=2^{-52}
\]
is the machine epsilon.

So far we have observed the following phenomena:

\begin{lstlisting}
>> a=2^53; b=(a-1)+1; isequal(a,b)
ans = 1
>> a=2^53; b=(a+1)-1; isequal(a,b)
ans = 0
>> a=1+2^-52; isequal(a,1)
ans = 0
>> a=1+2^-53; isequal(a,1)
ans = 1
\end{lstlisting}

Although matrix multiplication is associative over the real numbers, this property may fail in floating-point arithmetic because the underlying numerical operations are themselves affected by rounding, overflow, and underflow. For example,

\begin{verbatim}
a = 1e308;
b = 1e-308;
c = 1e-308;

(a*b)*c
a*(b*c)
\end{verbatim}
produces

\begin{verbatim}
1.0000e-308
0
\end{verbatim}
since \(b\,c=10^{-616}\) underflows to zero, whereas \(a\,b=1\) is represented exactly. Thus
\(
(a b)c \neq a(b c)
\)
in machine arithmetic. Since matrix multiplication is computed from floating-point additions and multiplications, associativity should likewise be regarded as an approximation rather than an exact identity in computational settings.

The philosophy adopted throughout this paper is that certain mathematical structures are best understood computationally through finite representations. Questions concerning machine precision and numerical stability, while important, lie outside the scope of the present work and will be addressed elsewhere. In this setting, sets are represented by row matrices, functions by indexing vectors, relations by logical matrices, and graphs by sparse matrices. The emphasis is therefore not on exact mathematical representations, but on computational models that preserve the essential categorical structure while remaining amenable to algorithmic manipulation.

Many of the examples will involve finite, dyadic, and combinatorial structures represented through finite numerical approximations. A typical example is the Boolean algebra \[ \{0,1\}^{\mathbb N}, \] whose elements may be viewed as infinite binary strings. Computationally, we approximate this space by the interval \[ [0.5,1), \] represented with double-precision floating-point numbers. Since IEEE 754 arithmetic allows us to control only the first fifty-two binary digits of a mantissa, each machine number may be interpreted as describing an infinite binary sequence whose initial segment is known while the remaining digits are left unspecified. This point of view allows one to retain a useful geometric picture. The interval provides a continuous model of an intrinsically discrete object, and dyadic rational numbers play a distinguished role. 

Following a well-known ambiguity of binary expansions, each dyadic number may be regarded as possessing two representatives, one terminating and one repeating. In some of our visualizations, dyadic points are therefore separated into conjugate pairs lying infinitesimally above and below the real axis, producing a geometric representation of the underlying Boolean structure while preserving its dyadic hierarchy (see Appendix \ref{Appendix Bollean} for further details). A similar philosophy will be adopted elsewhere.

 For example, finite-dimensional vector spaces may be represented computationally only up to the dimensions supported by the underlying machine arithmetic, while infinite-dimensional spaces such as \[ \mathbb R[x] \qquad\text{or}\qquad \mathbb R^{\mathbb N} \] are represented symbolically through finite descriptions of their finite dimensional approximations. The purpose of these constructions is not to reproduce infinite mathematics exactly, but rather to provide finite computational models that preserve enough structure to support categorical reasoning and algorithmic experimentation.

The reader is not expected to be an expert Octave programmer. The commands used throughout the paper are intentionally elementary and rely primarily on matrix operations, indexing, logical arrays, sparse matrices, and a small collection of built-in functions such as \texttt{unique}, \texttt{ismember}, \texttt{find}, \texttt{accumarray}, and \texttt{sparse}. One of the main themes of the paper is that these familiar operations already implement a surprisingly rich collection of categorical constructions. In Appendix \ref{sec: Appendix Octave} the reader will find a brief outline on basic programming with Matlab and Octave.

This observation naturally leads to the idea of a \emph{computational category}. Roughly speaking, such a category is one whose objects, morphisms, and categorical constructions admit concrete computational representations and can therefore be manipulated algorithmically. The precise definition will be developed later. For the moment, the reader may simply keep in mind the familiar examples of finite sets represented by indexing structures and finite-dimensional vector spaces represented by matrices.

\section{Examples of Computational Categories}

We conclude this introductory discussion with several computational categories that will serve as running examples throughout the paper. Each of them admits a description in terms of finite data structures and may therefore be regarded as a concrete realization of the general framework to be introduced later on.

\subsection*{The Category \(\mathbf{M1}\)}

The category \(\mathbf{M1}\) is the ordinary category of matrices.

Its objects are natural numbers

\[
0,1,2,\ldots,
\]
and a morphism

\[
A:n\to m
\]
is an \(m\times n\) matrix.

Composition is ordinary matrix multiplication,

\[
B\circ A=BA,
\]
and identities are the identity matrices

\[
I_n:n\to n.
\]

Thus \(\mathbf{M1}\) is precisely the representative category associated with finite-dimensional vector spaces equipped with ordered bases. 

In concrete implementations, however, the limitations of machine arithmetic must be taken into account. In particular, the dimensions appearing as objects of \(\mathbf{M1}\) are constrained by the numerical representation available in the underlying system. For MATLAB and Octave, which are based primarily on double-precision arithmetic, this means that dimensions can only be represented exactly up to \(2^{53}\). 
Similar considerations apply to other computational environments and should be regarded as implementation details rather than intrinsic limitations of the mathematical theory. Nevertheless, such restrictions become relevant when one interprets \(\mathbf{M1}\) as a concrete computational category. 

In future work, the consequences of imposing an upper bound such as \(2^{53}\) on object dimensions should be investigated more systematically. For the purposes of the present paper, we shall ignore this constraint and work formally with arbitrary finite dimensions. It is worth noting, however, that the restriction would genuinely alter the categorical properties of the resulting computational category. For example, if dimensions were required to satisfy \(n\leq 2^{53}\), then the product of two objects of dimension \(2^{32}\) would fail to exist, since \[ 2^{32}\cdot 2^{32}=2^{64}>2^{53}. \] Thus a finite computational realization need not inherit all limits and colimits of its mathematical counterpart.

\subsection*{The Category \(\mathbf{M2}\)}

The category \(\mathbf{M2}\) has the same objects as \(\mathbf{M1}\),

\[
0,1,2,\ldots,
\]
but a morphism

\[
A:n\to m
\]
is represented by an \(n\times m\) matrix. Note that in this category morphisms are displayed as
\[
\xymatrix{n\ar[r]^{A_{n\times m}} & m}
\]
whereas in $\mathbf{M1}$ the morphisms are displayed as
\[
\xymatrix{m\ar[r]^{A_{n\times m}} & n}
\]

Composition is defined by

\[
A\circ B=BA,
\]
so that matrix multiplication appears in the reverse order from the conventional one.

Equivalently,

\[ \mathbf{M2}\cong \mathbf{M1}^{op}. \] 

This category is often convenient when matrices are viewed as acting on row vectors rather than on column vectors. 

The transpose operation \[ A\longmapsto A' \] defines a contravariant equivalence \[ (-)' : \mathbf{M1}\longrightarrow\mathbf{M1}, \] or equivalently a covariant functor \[ (-)' : \mathbf{M1}\longrightarrow\mathbf{M2}, \] since \[ (BA)'=A'B'. \] In this sense, ordinary matrix transposition is the computational manifestation of categorical duality.

\subsection*{The Category \(\mathbf{Idx}\)}

The category \(\mathbf{Idx}\) is the computational category underlying finite sets.

Its objects are natural numbers

\[
0,1,2,\ldots,
\]
and a morphism

\[
u:n\to m
\]
is an indexing vector of length \(n\) whose entries lie in the range

\[
1,\ldots,m.
\]

Composition is ordinary indexing:

\[
(v\circ u)(i)=v(u(i)).
\]
In Octave notation this becomes

\begin{verbatim}
w = v(u);
\end{verbatim}
and identities are given by

\begin{verbatim}
(1:n)'
\end{verbatim}
for each object \(n\).

This category will play a central role in the sequel.

\subsection*{The Category \(\mathbf{Cols}\)}

The category \(\mathbf{Cols}\) is built from finite-dimensional column spaces.

Its objects are matrices \(A\) satisfying

\[
\operatorname{rank}(A)
=
\operatorname{size}(A,2).
\]

Equivalently, the columns of \(A\) are linearly independent and therefore form a basis of the corresponding column space.

A morphism

\[
T:A\to B
\]
is a matrix \(T\) of appropriate size such that there exists a matrix \(F\) satisfying

\[
BF=TA.
\]

Since the columns of \(B\) are linearly independent, the matrix \(F\) is uniquely determined and may be computed by

\begin{verbatim}
F = B\(T*A);
\end{verbatim}
The columns of each object already form a basis of the corresponding column space. Thus, unlike finite sets or finite-dimensional vector spaces, no auxiliary choice of ordering or basis is required. The objects carry their own coordinates.

From a computational perspective, \(\mathbf{Cols}\) may be viewed as an excellent Octave approximation to the category of finite-dimensional vector spaces. Every object is represented directly by a matrix, while morphisms are represented by matrices satisfying \(BF=TA\). Passing to the quotient by isomorphism identifies matrices that differ only by a change of coordinates, yielding a category equivalent to \(\mathbf{M1}\).

More explicitly, there are natural functors

\[
Q:\mathbf{Cols}\to\mathbf{M1}
\]
and

\[
I:\mathbf{M1}\to\mathbf{Cols}.
\]

For an object \(A\) of \(\mathbf{Cols}\), define

\[
Q(A)=\operatorname{size}(A,2),
\]
the number of columns of \(A\). Given a morphism

\[
T:A\to B,
\]
the defining relation

\[
BF=TA
\]
determines a unique matrix

\[
F=B\backslash(TA),
\]
which we take as the image of \(T\) under \(Q\). Thus

\[
Q(T)=F.
\]

On the other hand, every object \(n\) of \(\mathbf{M1}\) determines the matrix

\[
I(n)=I_n,
\]
and every morphism

\[
F:n\to m
\]
is sent to the morphism

\[
F:I_n\to I_m,
\]
since

\[
I_mF=FI_n.
\]

The composite

\[
Q\circ I
\]
is exactly the identity functor on \(\mathbf{M1}\). Moreover, for every object \(A\) of \(\mathbf{Cols}\), the columns of \(A\) define an isomorphism

\[
I(Q(A))
=
I_{\operatorname{size}(A,2)}
\longrightarrow A,
\]
and therefore \(I\circ Q\) is naturally isomorphic to the identity functor on \(\mathbf{Cols}\).

Consequently, \(\mathbf{M1}\) may be regarded as the representative category associated with \(\mathbf{Cols}\). Equivalently, \(\mathbf{Cols}\) is obtained from \(\mathbf{M1}\) by replacing the canonical basis \(I_n\) with arbitrary bases represented by full-rank matrices, while \(\mathbf{M1}\) is recovered by quotienting these choices of coordinates.

\subsection*{The Category \(\mathbf{Rows}\)}

The category \(\mathbf{Rows}\) provides a computational realization of finite sets.

Its objects are matrices whose rows are distinct and lexicographically sorted.

If \(A\) and \(B\) are such objects, a morphism

\[
v:A\to B
\]
is an indexing vector \(v\) of length

\[
\operatorname{size}(A,1)
\]
whose entries lie in the range

\[
1,\ldots,\operatorname{size}(B,1).
\]

The action of \(v\) is given by row selection: \[ A(i,:) \longmapsto B(v(i),:). \] In other words, the \(i\)-th element of the finite set represented by \(A\) is mapped to the \(v(i)\)-th element of the finite set represented by \(B\). Since elements are encoded by rows, this action is realized computationally by row indexing.

Composition is again ordinary indexing,

\[
w=u(v),
\]
and identities are given by

\begin{verbatim}
(1:size(A,1))'
\end{verbatim}
for each object \(A\).

The relationship between \(\mathbf{Rows}\) and \(\mathbf{Idx}\) is completely analogous to that between \(\mathbf{Cols}\) and \(\mathbf{M1}\). Indeed, there are functors

\[
Q:\mathbf{Rows}\to\mathbf{Idx}
\]
and

\[
I:\mathbf{Idx}\to\mathbf{Rows}.
\]

For an object \(A\), define

\[
Q(A)=\operatorname{size}(A,1),
\]
the number of rows of \(A\). A morphism

\[
v:A\to B
\]
is already an indexing vector, and therefore

\[
Q(v)=v.
\]

Conversely, for an object \(n\) of \(\mathbf{Idx}\), let

\[
I(n)=
\begin{bmatrix}
1\\
2\\
\vdots\\
n
\end{bmatrix},
\]
viewed as an object of \(\mathbf{Rows}\). Given a morphism

\[
v:n\to m,
\]
the functor \(I\) sends \(v\) to the same indexing vector regarded as a morphism

\[
I(n)\to I(m).
\]

The composite

\[
Q\circ I
\]
is therefore the identity functor on \(\mathbf{Idx}\). Moreover, every object \(A\) of \(\mathbf{Rows}\) is canonically isomorphic to \(I(Q(A))\), since its rows are already indexed by

\[
1,\ldots,\operatorname{size}(A,1).
\]

Consequently,

\[
I\circ Q
\]
is naturally isomorphic to the identity functor on \(\mathbf{Rows}\), and the two categories are equivalent:

\[
\mathbf{Rows}\simeq\mathbf{Idx}.
\]

Thus \(\mathbf{Rows}\) may be viewed as a concrete realization of finite sets in which elements are represented by rows of matrices and functions by indexing vectors, while \(\mathbf{Idx}\) provides the corresponding representative category whose objects are natural numbers.

The examples considered so far suggest a common pattern. Although they arise from different mathematical contexts, all of them admit concrete representations in which objects are encoded by finite data and morphisms by matrices, indexing vectors, or related computational structures. In each case, the original category is replaced by a collection of canonical representatives that can be manipulated directly within a computational environment. The purpose of the next section is to formalize this phenomenon.

\section{Computational Categories}

Many categories occurring in mathematics admit canonical representatives indexed by the natural numbers. Typical examples include finite sets, where the representative of dimension \(n\) is the set \(\{1,\ldots,n\}\), and finite-dimensional vector spaces, where the representative of dimension \(n\) is \(K^n\). The purpose of this section is to formalize the common structure underlying such examples.

Let \(\mathbf{C}\) be a category and let \[ \operatorname{emb}:\mathbf N\to \mathbf{C} \] be an embedding of the linearly ordered category of natural numbers $\mathbf{N}$ into the given category \(\mathbf{C}\). The objects \[ \operatorname{emb}(0),\operatorname{emb}(1),\operatorname{emb}(2),\ldots \] will be regarded as canonical representatives.

 An object \(X\in \mathbf{C}\) will be called \emph{finite} whenever there exists an integer \(n\) and an isomorphism

\[
\theta:\operatorname{emb}(n)\xrightarrow{\cong}X.
\]
Occasionally, we may also write $\operatorname{dim}(X)=n$.

The category \(\operatorname{Fin}(\mathbf{C})\) is defined as follows. Its objects are triples

\[
(X,n,\theta),
\qquad
\theta:\operatorname{emb}(n)\xrightarrow{\cong}X,
\]
and its morphisms

\[
(X,n,\theta)\longrightarrow(Y,m,\varphi)
\]
are simply the morphisms

\[
f:X\to Y
\]
of the original category \(\mathbf{C}\). Notice that no compatibility between \(f\) and the chosen isomorphisms \(\theta\) and \(\varphi\) is required. The coordinatizations are attached only to objects; morphisms remain exactly those of \(\mathbf{C}\). In a sense, the isomorphisms \(\theta\) and \(\varphi\) behave like hidden variables: they are present in the background and determine the coordinate description of the objects, but they do not appear explicitly in the morphisms. 

The chosen isomorphisms become relevant only when transporting morphisms to the canonical representatives. Indeed, every morphism

\[
f:(X,n,\theta)\to(Y,m,\varphi)
\]
determines a morphism

\[
\varphi^{-1}f\theta:
\operatorname{emb}(n)\to\operatorname{emb}(m).
\]

Let \(\sim\) denote the equivalence relation on objects of \(\operatorname{Fin}(\mathbf{C})\) defined by

\[
(X,n,\theta)\sim(Y,m,\varphi)
\]
whenever \(X\) and \(Y\) are isomorphic in \(\mathbf{C}\). The corresponding quotient category will be denoted by

\[
\operatorname{Rep}(\mathbf{C})
=
\operatorname{Fin}(\mathbf{C})/{\sim}.
\]

Its objects may therefore be identified with the natural numbers

\[
0,1,2,\ldots,
\]
corresponding to the chosen representatives \(\operatorname{emb}(n)\). Moreover, for every pair \(n,m\), the hom-set of \(\operatorname{Rep}(\mathbf{C})\) is given by

\[
\operatorname{Hom}_{\operatorname{Rep}(\mathbf{C})}(n,m)
=
\operatorname{Hom}_{\mathbf{C}}(\operatorname{emb}(n),\operatorname{emb}(m)).
\]

Thus \(\operatorname{Rep}(\mathbf{C})\) replaces finite objects by their dimensions while retaining the morphisms between the corresponding representatives.

The construction above captures finite objects only. In many situations, however, one wishes to keep track of a distinguished class of non-finite objects as well. Let \(V\) be an object of \(\mathbf{C}\) for which no isomorphism

\[
\operatorname{emb}(n)\xrightarrow{\cong}V
\]
exists. Let us say,

\[
\dim(V)=\mathrm{NaN}.
\]

Define \(\operatorname{Fin}(\mathbf{C},V)\) to be the category whose objects are either triples

\[
(X,n,\theta),
\qquad
\theta:\operatorname{emb}(n)\xrightarrow{\cong}X,
\]
or pairs

\[
(U,\varphi),
\qquad
\varphi:U\xrightarrow{\cong}V.
\]

Morphisms are those morphisms of \(\mathbf{C}\) satisfying

\[
\dim(\operatorname{dom}(f))=\mathrm{NaN}
\qquad\Longrightarrow\qquad
f \text{ is an isomorphism}.
\]

Thus non-finite objects contribute only through a chosen representative \(V\) and its automorphisms, while arbitrary morphisms from finite objects into \(V\) are still allowed.

Passing to the quotient by isomorphism classes yields a category

\[
\operatorname{Rep}(\mathbf{C},V)=\operatorname{Fin}(\mathbf{C},V)/{\sim},
\]
whose objects may be identified with

\[
0,1,2,\ldots,V.
\]

The hom-sets between finite objects remain

\[
\operatorname{Hom}_{\operatorname{Rep}(\mathbf{C},V)}(n,m)
=
\operatorname{Hom}_{\mathbf{C}}(\operatorname{emb}(n),\operatorname{emb}(m)),
\]
while the additional data are provided by

\[
\operatorname{Hom}_{\mathbf{C}}(\operatorname{emb}(n),V)
\]
and

\[
\operatorname{Aut}_{\mathbf{C}}(V).
\]

A category \(\mathbf{C}\) will be called \emph{computational of finite type} whenever

\[
\mathbf{C}/{\sim}
\cong
\operatorname{Rep}(\mathbf{C}),
\]
and \emph{computational of type \(V\)} whenever

\[
\mathbf{C}/{\sim}
\cong
\operatorname{Rep}(\mathbf{C},V).
\]

In the finite case, every isomorphism class is represented by a dimension

\[
0,1,2,\ldots.
\]

In the second case, one adjoins a distinguished representative \(V\) for a chosen class of non-finite objects.

For example, let \(\mathbf{C}=\mathbf{Vect}_{\mathbb R}\) and

\[
\operatorname{emb}(n)=\mathbb R^n.
\]

Then \(\operatorname{Rep}(\mathbf{C})\) describes finite-dimensional vector spaces through matrices and dimensions. Choosing

\[
V=\mathbb R[x]
\]
produces a computational category containing finite-dimensional vector spaces together with a representative for countably-dimensional vector spaces. Choosing instead

\[
V=\mathbb R^{\mathbb N}
\]
yields a different computational category in which the distinguished non-finite object has a larger dimension. Although one cannot generally represent all infinite dimensions simultaneously, many categories may be studied by selecting a suitable representative \(V\) and adjoining it to the finite part of the theory.

In both situations, objects are described by dimensions or by a distinguished representative \(V\), while morphisms are described by coordinate expressions relative to the representatives. Consequently, categorical constructions may be transported to \(\operatorname{Rep}(\mathbf{C})\) or \(\operatorname{Rep}(\mathbf{C},V)\), where they become amenable to algorithmic implementation.


The real strength of this viewpoint lies in its ability to combine finite and non-finite structures within a single computational setting. Rather than treating infinite objects as inaccessible, one selects representative objects through which they can interact with finite combinatorial data via explicit morphisms. This phenomenon is illustrated in Example~\ref{e.g. mixed morphism}, where finite sets, finite-dimensional vector spaces, and maps between them are brought together in the category \(\mathbf{Idx}+\infty\).

\begin{example}[Finite Sets]
Let \(\mathbf{Ens}_f\) denote the category of finite sets and define

\[
\operatorname{emb}(n)=\{1,\ldots,n\}.
\]

To specify an isomorphism

\[
\theta:\{1,\ldots,n\}\xrightarrow{\cong}X
\]
one must choose an ordering of the elements of \(X\). Thus \(\operatorname{Fin}(\mathbf{Ens}_f)\) consists of finite sets equipped with enumerations. Passing to \(\operatorname{Rep}(\mathbf{Ens}_f)\) identifies all enumerated sets of cardinality \(n\) with the representative \(\{1,\ldots,n\}\). Morphisms are then represented by indexing vectors. This is precisely the computational model underlying finite-set computations. In MATLAB and Octave, it is realized by the category \(\mathbf{Idx}\), whose objects are natural numbers and whose morphisms are indexing vectors composed through ordinary indexing.
\end{example}

\begin{example}[Finite-Dimensional Vector Spaces]
Let \(\mathbf{Vect}^{fd}_K\) denote the category of finite-dimensional vector spaces over a field \(K\) and define

\[
\operatorname{emb}(n)=K^n.
\]

Choosing an isomorphism

\[
\theta:K^n\xrightarrow{\cong}V
\]
is equivalent to choosing an ordered basis of \(V\). Consequently, \(\operatorname{Fin}(\mathbf{Vect}^{fd}_K)\) consists of finite-dimensional vector spaces equipped with ordered bases. After passing to \(\operatorname{Rep}(\mathbf{Vect}^{fd}_K)\), objects are identified with natural numbers and morphisms become matrices. Composition is represented by matrix multiplication. Thus \(\operatorname{Rep}(\mathbf{Vect}^{fd}_K)\) is precisely the category \(\mathbf{M1}\). In this sense, \(\mathbf{M1}\) constitutes the natural computational realization of finite-dimensional linear algebra and provides a particularly good model for MATLAB and Octave, whose basic data structures and operations are already matrix-oriented.
\end{example}

\begin{example}[Column Spaces]
Consider the category $\mathbf{Cols}$ from above.
The columns of an object \(A\) are linearly independent and therefore form a basis of its column space. Unlike the case of finite sets, where an ordering must be chosen, or finite-dimensional vector spaces, where a basis must be selected, no additional coordinatization is required here: the coordinates are already built into the object itself. Consequently, each object determines a canonical isomorphism \[ K^n \xrightarrow{\cong} \operatorname{Col}(A), \] where \(n\) is the number of columns of \(A\). The corresponding representative category is therefore equivalent to \(\mathbf{M1}\), whose objects are natural numbers and whose morphisms are matrices. Thus this category is computational in a particularly strong sense: not only are its objects and morphisms completely described by matrices, but the objects already carry their own coordinates.
\end{example}

\begin{example}[Finite Sets and Vector Spaces]\label{e.g. mixed morphism}

The previous examples may be combined into a single computational category containing both finite sets and canonical finite-dimensional vector spaces.

The motivation comes from applications where one must manipulate, within the same framework, finite subsets of vector spaces together with maps between them, while also allowing maps from finite subsets into the set of all vectors of a vector space. In such situations, one encounters three different kinds of morphisms:

\begin{enumerate}
\item maps between finite subsets;
\item maps from finite subsets into an ambient vector space;
\item isomorphisms between ambient vector spaces.
\end{enumerate}

To model this situation, consider the category whose objects are subsets of vector spaces, including the vector spaces themselves. As representatives, we choose the finite sets

\[
\{1,\ldots,n\}
\]
and the canonical vector spaces

\[
K^m.
\]

This leads naturally to two kinds of objects, represented respectively by

\[
0,1,2,\ldots
\]
and

\[
\infty+0,\infty+1,\infty+2,\ldots.
\]

The rule that every morphism having a non-finite domain must be an isomorphism introduces an asymmetry: morphisms from finite objects into vector spaces are allowed, whereas morphisms in the opposite direction are not. Moreover, the only morphisms between non-finite objects are isomorphisms. For simplicity, these isomorphisms will be required to be linear and therefore will be represented by invertible matrices.

The resulting category will be denoted by

\[
\mathbf{Idx}+\infty.
\]

Its objects are of two kinds.

Finite objects are represented by the natural numbers

\[
0,1,2,\ldots,
\]
and should be interpreted as the cardinalities of finite subsets of vector spaces, encoded computationally by row matrices.

For each \(m\geq0\), the object

\[
\infty+m
\]
represents the canonical vector space

\[
K^m.
\]

Morphisms between finite objects

\[
n\to p
\]
are indexing vectors, exactly as in \(\mathbf{Idx}\).

Morphisms between non-finite objects

\[
\infty+n
\longrightarrow
\infty+n
\]
are invertible matrices. Thus

\[
\operatorname{Aut}(\infty+n)
=
GL_n(K).
\]

Finally, mixed morphisms

\[
n
\longrightarrow
\infty+m
\]
are represented by row matrices

\[
F=
\begin{bmatrix}
f_1\\
\vdots\\
f_n
\end{bmatrix},
\]
having \(m\) columns. Equivalently, such a morphism assigns to each element of the finite set \(n\) a vector of \(K^m\). The rows of \(F\) are therefore precisely the images of the elements

\[
1,\ldots,n.
\]

Composition with an indexing vector

\[
u:p\to n
\]
is given by row selection,

\[
F\circ u = F(u,:),
\]
while composition with an automorphism

\[
A:\infty+m\to\infty+m
\]
is given by

\[
A\circ F = FA'.
\]

Thus finite combinatorial data and linear-algebraic data coexist within the same category. Finite sets are represented by indexing structures, vector spaces by matrices, and mixed morphisms by finite collections of vectors. As we shall see later, this construction provides a convenient computational model for relating finite geometric data to linear-algebraic structures.

\end{example}

Note that, in an Octave implementation, the object \(\infty+m\) may be encoded by the complex number \texttt{Inf+1i*m}.

\section{A Quaternionic Example}

The previous examples were based on finite sets and finite-dimensional vector spaces. We now consider a computational category whose objects are subsets of the quaternion algebra

\[
\mathbb H.
\]

Let \(\mathbf{C}\) denote the category whose objects are subsets

\[
A\subseteq\mathbb H,
\]
and whose morphisms are ordinary functions between subsets.

As before, we choose an embedding

\[
\operatorname{emb}:\mathbf N\to \mathbf{C},
\]
by setting

\[
\operatorname{emb}(n)=\{1,\ldots,n\},
\]
viewed as a subset of \(\mathbb H\). Finite objects are therefore finite subsets of \(\mathbb H\).

An object of \(\operatorname{Fin}(\mathbf{C})\) is a triple

\[
(A,n,\theta),
\qquad
\theta:\{1,\ldots,n\}\xrightarrow{\cong}A.
\]

Equivalently, it is a finite subset of \(\mathbb H\) equipped with an ordering of its elements. Such objects admit a natural matrix representation: arranging the coordinates of the quaternions into rows yields a matrix

\[
A=
\begin{bmatrix}
a_1 & b_1 & c_1 & d_1\\
a_2 & b_2 & c_2 & d_2\\
\vdots & \vdots & \vdots & \vdots\\
a_n & b_n & c_n & d_n
\end{bmatrix},
\]
whose rows are distinct and are assumed to be sorted lexicographically.

Thus \(\operatorname{Fin}(\mathbf{C})\) may be identified with a subcategory of
\(\mathbf{Rows}\) introduced earlier. Objects are matrices with distinct
sorted rows (and exactly four columns) and morphisms are indexing vectors

\[
v:\{1,\ldots,n\}\to\{1,\ldots,m\},
\]
implemented computationally by row selection.

Passing to the quotient by isomorphism classes yields

\[
\operatorname{Rep}(\mathbf{C}),
\]
whose objects are the natural numbers

\[
0,1,2,\ldots,
\]
and whose morphisms are indexing vectors. Consequently,

\[
\operatorname{Rep}(\mathbf{C})\cong\mathbf{Idx}.
\]

We now extend the construction by choosing

\[
V=\mathbb H.
\]

Since \(\mathbb H\) is not finite, we have

\[
\dim(\mathbb H)=\mathrm{NaN}.
\]

The category \(\operatorname{Fin}(\mathbf{C},\mathbb H)\) contains not only finite
ordered subsets of \(\mathbb H\), but also objects

\[
(U,\varphi),
\qquad
\varphi:U\xrightarrow{\cong}\mathbb H.
\]

Up to isomorphism, these all represent the single object \(\mathbb H\).

The corresponding representative category

\[
\operatorname{Rep}(\mathbf{C},\mathbb H)
\]
therefore has objects

\[
0,1,2,\ldots,\mathbb H.
\]

The morphisms between finite objects are again indexing vectors.
A morphism

\[
n\to\mathbb H
\]
is represented by a row matrix

\[
A=
\begin{bmatrix}
a_1 & b_1 & c_1 & d_1\\
a_2 & b_2 & c_2 & d_2\\
\vdots & \vdots & \vdots & \vdots\\
a_n & b_n & c_n & d_n
\end{bmatrix},
\]
whose rows describe the images of the elements

\[
1,\ldots,n
\]
in \(\mathbb H\).

Unlike the objects of \(\mathbf{Rows}\), these rows are not required to be
distinct or sorted. Indeed, arbitrary functions

\[
\{1,\ldots,n\}\to\mathbb H
\]
are allowed, and repeated rows correspond precisely to repeated values of
the function.

The remaining morphisms are provided by the automorphism group

\[
\operatorname{Aut}_{\mathbf{C}}(\mathbb H).
\]

The structure of this group will not be described here, although it is
present in the representative category and acts on the morphisms
\[
n\to\mathbb H.
\]

Thus \(\operatorname{Rep}(\mathbf{C},\mathbb H)\) provides a computational model
consisting of finite dimensions together with a distinguished infinite
object. Finite subsets are handled through the category \(\mathbf{Rows}\),
while maps into the whole quaternion algebra are encoded directly by row
matrices. This example illustrates how a representative category may
combine finite combinatorial data with a chosen non-finite object within
the same computational framework.

\section{Examples of Application}

We now turn to a few simple applications illustrating how the computational categories introduced in the previous sections interact in practice.

As observed in \cite{MFACS2023,MFMath2023}, there are remarkable parallels between internal groupoid structures and triangulated surfaces. Although the details lie beyond the scope of the present paper, both may be viewed as arising from a common underlying structure, namely that of a \emph{multilink}. Roughly speaking, a multilink consists of an indexing set \(A\), a collection of endomorphisms

\[
\varphi_i:A\to A,
\]
a family of projection maps

\[
s_j:A\to B,
\]
and a geometric realization

\[
g:A\to K^n,
\]
where \(K^n\) may represent Euclidean space, the complex numbers, the quaternions, or a more general geometric algebra. Such structures provide a common framework for studying both geometric and combinatorial constructions.

From this perspective, it becomes natural to consider simultaneously indexing maps

\[
v:n\to m
\]
and geometric maps

\[
F:n\to\mathbb H.
\]

The categories \(\mathbf{M1}\), \(\mathbf{M2}\), \(\mathbf{Cols}\), \(\mathbf{Rows}\), \(\mathbf{Idx}\), and the representative categories introduced earlier may therefore be regarded as different computational languages describing various aspects of the same underlying structures.

In the remainder of this section we examine several elementary constructions arising from their interaction. The reader may find it useful to consult the appendices for additional background and implementations.

\subsection{Image Factorizations}

One particularly important example is the image factorization. Given a morphism

\[
f:A\to B
\]
in \(\mathbf{Rows}\), represented by an indexing vector \(f\), the image of \(f\) may be computed through

\begin{verbatim}
[Q,r,q] = unique(B(f,:),'rows');
\end{verbatim}

The matrix \(Q\) contains the distinct rows occurring in the image of \(f\), the vector

\[
r:Q\to A
\]
selects one representative for each fibre of \(q\), and

\[
q:A\to Q
\]
is the canonical surjection. Thus \(f\) factors as

\[
A
\xrightarrow{q}
Q
\xrightarrow{f(r)}
B.
\]

The map \(q\) is surjective, while \(f(r)\) is injective.

The same command admits a second interpretation. Consider a morphism

\[
F:n\to\mathbb H
\]

represented by a row matrix

\[
F=
\begin{bmatrix}
a_1&b_1&c_1&d_1\\
\vdots&\vdots&\vdots&\vdots\\
a_n&b_n&c_n&d_n
\end{bmatrix}.
\]

Applying

\begin{verbatim}
[V,r,q] = unique(F,'rows');
\end{verbatim}
produces the factorization

\[
n
\xrightarrow{q}
\operatorname{size}(V,1)
\xrightarrow{V}
\mathbb H,
\]
where \(V\) consists of the distinct quaternionic values assumed by \(F\).

Thus \texttt{unique} may be interpreted either as the image factorization of a morphism in \(\mathbf{Rows}\) or as the image factorization of a geometric map into \(\mathbb H\). In both cases, the command separates the essential geometric data from the indexing information required to reconstruct the original map. This dual interpretation is characteristic of computational categories: the same algorithm simultaneously carries geometric and categorical meaning.

\subsection{Graphs and relations}

A second family of examples arises from relations. Let \(\mathbf{Rel}\) denote the category whose objects are logical matrices. Given relations \[ R:n\to m, \qquad S:m\to p, \] there is a natural horizontal composition, namely ordinary relational composition, which in matrix form is computed by \[ S\star R = \operatorname{logical}(R*S). \] This operation is fundamental and will be used repeatedly throughout the paper. However, it should not be confused with the actual composition of morphisms in the category \(\mathbf{Rel}\). Indeed, a morphism \[ (f,g):R\to S \] consists of morphisms \[ f:n_R\to n_S, \qquad g:m_R\to m_S \] in \(\mathbf{Idx}\) satisfying \[ R(i,j)=1 \qquad\Longrightarrow\qquad S(f(i),g(j))=1. \] In other words, every related pair in \(R\) is sent to a related pair in \(S\). The composition of such morphisms is simply inherited from \(\mathbf{Idx}\): \[ (f',g')\circ(f,g) = (f'\circ f,\,g'\circ g). \] Thus relations carry two natural compositions. One is the horizontal composition of relations themselves, given by Boolean matrix multiplication. The other is the ordinary categorical composition of structure-preserving morphisms between relations. This distinction is one of the first indications of the double-categorical nature of relations.

Relations may equally well be represented through the category \(\mathbf{Span}\). Given a logical matrix \(R\), the command \begin{verbatim} [d,c] = find(R); \end{verbatim} returns two indexing vectors having the same length, \[ d,c:\operatorname{nnz}(R)\to n,m, \] where \(\operatorname{nnz}(R)\) denotes the number of nonzero entries of \(R\). Computationally, the pair \[ (d,c) \] encodes the relation \(R\) by recording the row and column positions of its nonzero entries. Conversely, every pair of indexing vectors \[ d:E\to n, \qquad c:E\to m \] having common domain \(E\) determines a relation through \begin{verbatim} R = sparse(d,c,true,n,m); \end{verbatim} whose nonzero entries occur precisely at the positions specified by \(d\) and \(c\). This leads naturally to the category \(\mathbf{Span}\), whose objects are pairs \[ (d,c) \] of morphisms in \(\mathbf{Idx}\) having the same domain. Thus an object of \(\mathbf{Span}\) may be represented diagrammatically as \[ n \overset{d}{\longleftarrow} E \overset{c}{\longrightarrow} m, \] or computationally by the pair of indexing vectors \((d,c)\). The commands \begin{verbatim} [d,c] = find(R) \end{verbatim} and \begin{verbatim} R = sparse(d,c,true,n,m) \end{verbatim} therefore define functors \[ \mathbf{Rel} \;\xrightarrow{\;\texttt{find}\;} \mathbf{Span} \] and \[ \mathbf{Span} \;\xrightarrow{\;\texttt{sparse}\;} \mathbf{Rel}, \] which identify logical matrices with their corresponding span representations. In this sense, \(\mathbf{Rel}\) and \(\mathbf{Span}\) provide two computational descriptions of the same underlying data: logical matrices on one side and pairs of indexing vectors on the other.

From a practical point of view, the passage between
\(\mathbf{Span}\) and \(\mathbf{Rel}\) becomes particularly valuable
when one wishes to perform closure and quotient constructions.

Given a relation represented by a logical matrix \(R\), its
transitive closure may be computed either by a Warshall-style
algorithm or by iterative squaring. The latter is particularly
natural from the present perspective since it realizes the
transitive closure as the least fixed point of

\[
T\longmapsto T\cup(T\circ T).
\]

Starting with

\[
T=R,
\]
one iterates

\begin{verbatim}
Tnew = T | logical(T*T);
\end{verbatim}
until a fixed point is reached. When reflexivity has already
been imposed, one has

\[
T\subseteq T^2,
\]
and therefore the update simplifies to

\begin{verbatim}
Tnew = logical(T*T);
\end{verbatim}

This observation leads to substantial computational
improvements for reflexive-transitive and equivalence
closures.

An important application arises from coequalizers, as illustrated in  Appendix \ref{subsec: coeq}. Briefly, consider a
parallel pair

\[
d,c:E\rightrightarrows n
\]
in \(\mathbf{Idx}\). Applying \texttt{sparse} produces the
relation

\begin{verbatim}
R = sparse(c,d,true,n,n);
\end{verbatim}
where the order of \(c\) and \(d\) is chosen so that the pair

\[
(d(e),c(e))
\]
is interpreted as an identification to be imposed on the
quotient.

The smallest equivalence relation containing these pairs is
then obtained by computing the reflexive-symmetric-transitive
closure

\[
T=(R\cup R^{-1}\cup I)^*.
\]

The corresponding quotient is recovered through

\begin{verbatim}
[p1,p2] = find(triu(T));
t(p1) = p2;
\end{verbatim}

Since \texttt{find} scans the matrix in sorted order, the
vector \(t\) associates to each element the chosen
representative of its equivalence class. Thus \(t\) is
precisely the split coequalizer of the pair \((d,c)\), assuming that $d$ is surjective. See Appendix \ref{subsec: coeq} for further details.

 This illustrates another recurring theme of the paper: standard Octave commands such as \texttt{sparse}, \texttt{find}, and Boolean matrix multiplication naturally realize categorical constructions that would otherwise be described abstractly in terms of universal
properties.

\section{The Graph Construction}

One of the most useful applications of the categories introduced in the previous sections is the conversion of geometric data into graph-theoretic structures. The purpose of the \texttt{graph} construction is precisely to transform a pair of maps

\[
\texttt{source},\texttt{target}:n_E\to\mathbb H
\]
represented by row matrices into an internal graph in the category \(\mathbf{Idx}\).

The starting point is a collection of edges described geometrically. Thus the rows

\[
\texttt{source}(i,:)
\qquad\text{and}\qquad
\texttt{target}(i,:)
\]
represent the source and target of the \(i\)-th edge. No indexing of vertices is initially assumed; the graph is specified entirely through points of \(\mathbb H\).

For example,

\[
\texttt{source}=
\begin{bmatrix}
1&0&0&0\\
0&1&0&0\\
0&0&1&0
\end{bmatrix},
\qquad
\texttt{target}=
\begin{bmatrix}
0&1&0&0\\
0&0&1&0\\
1&0&0&0
\end{bmatrix}
\]
describe a directed cycle on three quaternionic vertices.

The first step consists in replacing geometric vertices by indexed ones. This is achieved through the image factorization

\begin{verbatim}
[V,e_,d] = unique(source,'rows');
\end{verbatim}

The matrix \(V\) contains the distinct source vertices, while

\[
d:n_E\to n_V
\]
assigns to each edge the index of its source vertex. Thus

\[
\texttt{source}
=
V(d,:).
\]

The command \texttt{unique} therefore constructs a vertex set together with a canonical indexing.

The second step identifies those target vertices already represented in \(V\):

\begin{verbatim}
[k,c] = ismember(target,V,'rows');
\end{verbatim}

For every edge satisfying \(k(i)=1\),

\[
\texttt{target}(i,:)
=
V(c(i),:).
\]

The vector \(c\) is therefore a partial target map. Since
\texttt{ismember} returns \(0\) whenever a target row does not occur in
\(V\), one restricts attention to the subset

\[
E'
=
\{\,e\in n_E\mid k(e)=1\,\}
\]
of edges whose targets are defined.

The graph itself is then obtained by

\begin{verbatim}
s = d(k);
t = c(k);
\end{verbatim}
producing a pair of morphisms

\[
s,t:n_{E'}\to n_V.
\]

In other words, the geometric data have been transformed into an internal graph

\[
n_{E'}
\overset{s}{\underset{t}{\rightrightarrows}}
n_V
\]
in the category \(\mathbf{Idx}\).

Observe that this graph is no longer described by points of
\(\mathbb H\). All geometric information has been replaced by indexing
data. The passage

\[
(\texttt{source},\texttt{target})
\longmapsto
(s,t)
\]
is therefore a conversion from coordinates to morphisms.

As can be seen from the implementation

\begin{lstlisting}
[V,e_,d]=unique(source,'rows');
[k,c]=ismember(target,V,'rows');
s=d(k);
t=c(k);
Adj=sparse(s,t,E(k));
cval=zeros(size(k)); cval(k)=1;
cval(~k)=-1; cval(e_)=2;
Edges=sparse(d,E,cval);
\end{lstlisting}
the first step computes the image factorization of the source map. Indeed,

\[
n_E
\xrightarrow{d}
n_V
\xrightarrow{V}
\mathbb H
\]
is the factorization produced by \texttt{unique}, where \(d\) is a split epimorphism, \(V\) is a monomorphism, and

\[
\texttt{source}=V(d,:).
\]

The second step computes the pullback of the target map along the monomorphism \(V\). Since \(V\) is injective, every target vertex belonging to the image of \(V\) factors uniquely through \(V\). The command

\begin{verbatim}
[k,c]=ismember(target,V,'rows');
\end{verbatim}
records precisely this information. The logical vector \(k\) determines the pullback object

\[
E'
=
\{\,e\in E\mid k(e)=1\,\},
\]
while the corresponding projections are encoded by

\[
\pi_1=\texttt{find}(k),
\qquad
\pi_2=c(k).
\]

The graph morphisms are then obtained simply by composing these projections with the source and target assignments:

\[
s=d(\texttt{find}(k)),
\qquad
t=\pi_2.
\]

Thus the graph construction is obtained by combining an image factorization with a pullback. Geometric data in \(\mathbb H\) are first replaced by indexed representatives, and the resulting graph is expressed entirely in terms of morphisms of \(\mathbf{Idx}\).

It is worth observing that the information discarded in the passage from geometry to graphs is not actually lost. The vector \texttt{e\_}, returned by \texttt{unique}, records the representatives chosen for the image factorization, while the original edge indices \(E\) are retained through the auxiliary matrix \texttt{Edges}. Consequently, the implementation keeps track not only of the graph itself but also of the process by which it was obtained.

Once the graph has been indexed, standard matrix representations become available. In the implementation this is achieved through \begin{verbatim} 
Adj = sparse(s,t,E(k)); 
\end{verbatim}
 which constructs an adjacency matrix whose nonzero entries retain the labels of the original edges. Thus the same graph admits three complementary descriptions: \[ (\texttt{source},\texttt{target}), \] \[ n_{E'} \overset{s}{\underset{t}{\rightrightarrows}} n_V, \] and \[ Adj. \] The implementation also produces the auxiliary matrix \[ Edges=sparse(d,E,cval), \] where the values stored in \texttt{cval} distinguish between edges that belong to the indexed graph, edges discarded because their targets do not appear in the chosen vertex set, and representative edges selected during the image factorization. Consequently, the construction records not only the resulting graph but also part of the history by which the graph was extracted from the original geometric data. The significance of the construction lies not in graph theory itself, but in the fact that it combines several categorical mechanisms already encountered in earlier sections: \begin{itemize} \item \texttt{unique} computes the image factorization of the source map; \item \texttt{ismember} computes the pullback of the target map along the monomorphism determined by the image factorization; \item \texttt{sparse} converts the resulting pair of morphisms into a relation or adjacency matrix. \end{itemize} The entire procedure may therefore be regarded as a functorial passage from geometric objects in \(\mathbb H\) to internal graphs in \(\mathbf{Idx}\). Rather than constructing a graph from scratch, the algorithm reveals categorical structure that was already implicit in the original geometric data while simultaneously retaining enough information to reconstruct how the graph was obtained.

\section{Conclusion} 

This paper introduced a collection of computational categories arising naturally from MATLAB and Octave data structures. Matrices, indexing vectors, logical matrices, sparse matrices, finite point clouds, and graph representations were shown to admit categorical interpretations in which objects and morphisms become directly computable. 

The central idea is that many categories arising in scientific computing admit representative forms whose objects are indexed by dimensions and whose morphisms are represented by finite data structures.

 This viewpoint led to the notion of a computational category and to several concrete examples, including matrix categories, index categories, categories of row and column representations, categories of relations and spans, and categories combining finite sets with finite-dimensional vector spaces through mixed morphisms. 
 
 Beyond the individual examples, the paper illustrates a broader principle: standard computational operations such as indexing, row selection, factorization through \texttt{unique}, pullback constructions through \texttt{ismember}, and relational constructions through \texttt{sparse} already implement a surprisingly rich collection of categorical constructions. Consequently, category theory may serve not only as a language for describing computation, but also as a practical tool for organizing and understanding computational structures. 
 
 Future work will focus on extending the present framework into a general theory of computational categories. Besides studying additional examples and implementation-dependent phenomena such as bounded dimensions and finite precision, a central objective is the development of a comprehensive software toolbox for computational category theory. Such a system would provide a common interface for representing computational categories, constructing their representative categories, and computing categorical constructions such as limits, colimits, factorization systems, relations, graphs, quotients, and adjunction-like structures. Ultimately, the aim is to bring categorical programming to the same level of accessibility that matrix-oriented programming currently enjoys in MATLAB and Octave.

\section*{Acknowledgements and References} \addcontentsline{toc}{section}{Acknowledgements and References} 


The author would like to thank the participants of ITÁCA Fest 2025 and ICMS 2026 (Waterloo, Canada), together with colleagues, students, and collaborators at the Polytechnic of Leiria, for many fruitful discussions that contributed to the development of the ideas presented in this paper. 

Special thanks are due to Stellenbosch University for the generous hospitality during August 2026, where two talks related to this work were delivered. The questions and discussions that followed played an important role in refining several of the concepts developed herein.

\appendix

\section{Computational Background}\label{sec: Appendix Octave}

The purpose of this appendix is twofold. On the one hand, it provides readers with enough familiarity with Octave to follow the computational examples used throughout the paper. On the other hand, it serves as a brief introduction to Octave for students encountering the language for the first time.

Octave is a high-level language designed for scientific computing and numerical mathematics. Its syntax is largely compatible with MATLAB, and both systems are fundamentally matrix-oriented. Unlike many programming languages, where arrays are merely one data structure among many, matrices constitute the primary objects of computation. Consequently, understanding matrices and indexing operations is the key to understanding the language itself.

\subsection{The Workspace}

An Octave session takes place within a \emph{workspace}, namely the collection of currently defined variables. For example,

\begin{verbatim}
a = 2;
b = 3;
c = a + b;
\end{verbatim}
creates three variables.

The commands

\begin{verbatim}
who
whos
\end{verbatim}
display the contents of the workspace and information about the stored variables. Each variable has a name, a class (type), dimensions, and a value. Conceptually, the workspace may be viewed as a finite collection of labels associated with mathematical objects.

\subsection{Vectors, Matrices, and Dimensions}

Vectors and matrices are the fundamental data structures in Octave.

A row vector is entered as

\begin{verbatim}
v = [1 2 3 4];
\end{verbatim}
while a column vector is entered as

\begin{verbatim}
v = [1;2;3;4];
\end{verbatim}
or 

\begin{verbatim}
v = [1 2 3 4]';
\end{verbatim}

Matrices are formed by separating rows with semicolons:

\begin{verbatim}
A = [1 2;
     3 4;
     5 6];
\end{verbatim}

The commands

\begin{verbatim}
size(A)
rows(A)
columns(A)
\end{verbatim}
return information about the dimensions of \(A\). Since many categorical constructions in the paper are expressed in terms of finite sets represented by matrices, these commands will appear frequently.

\subsection{Indexing}

The most important operation in Octave is indexing.

Given

\begin{verbatim}
v = [10 20 30 40];
\end{verbatim}
one obtains the second element through

\begin{verbatim}
v(2)
\end{verbatim}
while, for a matrix,

\begin{verbatim}
A(2,:)
A(:,1)
\end{verbatim}
return respectively the second row and the first column.

The symbol

\begin{verbatim}
:
\end{verbatim}
means ``all indices in that dimension''.

Several rows may also be selected simultaneously. For example,

\begin{verbatim}
u = [3;1;1;2];
A(u,:)
\end{verbatim}
returns rows \(3,1,1,2\) of \(A\).

Throughout the paper this operation will be interpreted categorically. Indexing vectors behave exactly like functions between finite sets.

\subsection{Indexing Vectors as Functions}

Consider

\begin{verbatim}
u = [2;1;3;3];
v = [3;1;2];
\end{verbatim}

Then

\begin{verbatim}
v(u)
\end{verbatim}
returns

\begin{verbatim}
[1;3;2;2]
\end{verbatim}
obtained by replacing each entry \(u(i)\) with \(v(u(i))\).

We may interpret

\[
u:(1:4)\to(1:3)
\]
and

\[
v:(1:3)\to(1:3)
\]
as functions between finite sets. Under this interpretation

\[
v(u)
=
v\circ u.
\]

This simple observation is the foundation of most categorical constructions developed in the main text.

\subsection{Logical Vectors and Masks}

Logical vectors represent subsets and filtering operations.

For example,

\begin{verbatim}
v = [10 20 30 40];
k = [true false true false];
\end{verbatim}
gives

\begin{verbatim}
v(k)
\end{verbatim}
with output

\begin{verbatim}
[10 30]
\end{verbatim}

The vector \(k\) is called a \emph{mask}. It selects those positions for which the corresponding logical value is true. Conceptually, masks behave like characteristic functions of subsets and occur frequently in commands such as \texttt{find} and \texttt{ismember}.

\subsection{Derived Variables}

Variables are often obtained from other variables through indexing.

For example,

\begin{verbatim}
A = [1 2;
     3 4;
     5 6];

B = A([3 1],:);
\end{verbatim}
constructs \(B\) by selecting rows \(3\) and \(1\) of \(A\).

From a mathematical perspective, only the resulting matrix \(B\) matters. From a computational perspective, the indexing vector records part of the history of the construction. The categorical viewpoint adopted in the paper focuses primarily on the resulting objects and morphisms.

\subsection{Concatenation}

Matrices may be combined by concatenation.

Horizontal concatenation:

\begin{verbatim}
[A B]
\end{verbatim}
places matrices side by side.

Vertical concatenation:

\begin{verbatim}
[A;
 B]
\end{verbatim}
places one matrix above another.

For example,

\begin{verbatim}
A = [1 2];
B = [3 4];

[A B]
\end{verbatim}
returns

\[
[1\;2\;3\;4].
\]

Such constructions later appear as computational realizations of categorical coproducts.

\subsection{Important Built-In Functions}

Several built-in commands play a central role throughout the paper.

\paragraph{\texttt{unique}}

\begin{verbatim}
[B,r,q] = unique(A,'rows');
\end{verbatim}
extracts the distinct rows of \(A\) and records how every row of \(A\) corresponds to a row of \(B\). This operation will later be interpreted as an image factorization.

\paragraph{\texttt{ismember}}

\begin{verbatim}
[k,p] = ismember(A,B,'rows');
\end{verbatim}
determines which rows of \(A\) occur among the rows of \(B\). It plays a central role in inverse-image and pullback constructions.

\paragraph{\texttt{find}}

\begin{verbatim}
[i,j] = find(R);
\end{verbatim}
returns the positions of the nonzero entries of a matrix. In graph-theoretic language, it converts adjacency matrices into edge lists.

\paragraph{\texttt{sparse}}

\begin{verbatim}
R = sparse(i,j,1);
\end{verbatim}
constructs a sparse matrix whose nonzero entries occur at positions \((i,j)\). This command provides one of the principal bridges between matrices, relations, and graphs.

\paragraph{\texttt{sortrows}}

\begin{verbatim}
sortrows(A)
\end{verbatim}
sorts the rows of a matrix lexicographically and is frequently used when constructing canonical representatives.

\paragraph{\texttt{accumarray}}

\begin{verbatim}
accumarray(i,v)
\end{verbatim}
aggregates values according to an indexing vector \(i\). It may be interpreted as collecting information along the fibres of a function.

\subsection{Sparse Matrices}

Ordinary matrices store every entry explicitly. Sparse matrices store only the nonzero entries.

For example,

\begin{verbatim}
R = full(sparse([1;3],[2;1],[1;1],3,3));
\end{verbatim}
creates

\[
\begin{bmatrix}
0&1&0\\
0&0&0\\
1&0&0
\end{bmatrix}.
\]

Sparse matrices are particularly important because many relations and graphs contain far fewer nonzero entries than the total number of matrix positions. In the categorical language developed throughout the paper, they provide a natural computational representation of finite relations and directed graphs.

\section{Limits, Colimits, and Factorizations in \(\mathbf{Idx}\)}

One of the advantages of the category \(\mathbf{Idx}\) is that many fundamental categorical constructions admit remarkably simple implementations. Since objects are natural numbers and morphisms are indexing vectors, limits, colimits, and factorizations may often be expressed directly through standard Octave commands.

\subsection*{Subobjects and quotients}

A subobject of \(n_A\) is represented by a monomorphism

\[
m:n_S\rightarrowtail n_A.
\]

In \(\mathbf{Idx}\), monomorphisms are precisely injective indexing vectors. In the particular case of a monotone monomorphism, a subobject may be represented equivalently by a logical mask

\[
k:n_A\to\{0,1\}.
\]

The corresponding inclusion is obtained through

\begin{verbatim}
m = find(k);
\end{verbatim}
which extracts precisely those positions at which \(k\) is true.

Conversely,

\begin{verbatim}
[k,p] = ismember(f,m);
\end{verbatim}
tests whether a morphism \(f\) factors through the subobject represented by \(m\). Thus \texttt{ismember} may be viewed as implementing inverse-image constructions and pullbacks along monomorphisms.

It is worth observing that a logical vector

\[
k:n_A\to\{0,1\}
\]
admits two complementary interpretations.

On the one hand, it determines a subobject through

\begin{verbatim}
m = find(k);
\end{verbatim}
which acts as a filter, selecting those elements for which \(k=1\). In this interpretation, the logical vector describes an inclusion

\[
n_S\hookrightarrow n_A.
\]

On the other hand, the same logical vector may be interpreted as a monotone surjection through

\begin{verbatim}
q = cumsum([1; k(:)]);
\end{verbatim}
or, after a suitable normalization,

\begin{verbatim}
q = cumsum(k);
\end{verbatim}

The effect is entirely different. Whenever \(k(i)=0\), the value of \(q\) remains unchanged and the corresponding element is identified with the previous one. Whenever \(k(i)=1\), a new value is created. Consequently, \(q\) defines a quotient object by collapsing neighbouring elements into equivalence classes.

Thus the logical vector may be viewed either as a filtering device or as a broadcasting device. The command

\begin{verbatim}
find
\end{verbatim}
extracts the selected elements and therefore describes a subobject, whereas

\begin{verbatim}
cumsum
\end{verbatim}
groups elements together and therefore describes a quotient object.

From a categorical perspective, these two constructions provide computational realizations of dual notions. The first produces a monomorphism and hence a subobject, while the second produces a canonical surjection and hence a quotient. In this sense, logical vectors form a simple bridge between limits and colimits in \(\mathbf{Idx}\).

\subsection*{Products}

The product of two objects

\[
n_1
\qquad\text{and}\qquad
n_2
\]
is simply

\[
n_1n_2.
\]

The projection morphisms are obtained through

\begin{verbatim}
[pi1,pi2] = ind2sub([n1,n2],(1:n1*n2)');
\end{verbatim}
yielding

\[
\pi_1:n_1n_2\to n_1,
\qquad
\pi_2:n_1n_2\to n_2.
\]

Given morphisms

\[
f:Z\to n_1,
\qquad
g:Z\to n_2,
\]
the unique pairing morphism

\[
\langle f,g\rangle:Z\to n_1n_2
\]
is computed by

\begin{verbatim}
h = sub2ind([n1,n2],f,g);
\end{verbatim}
and satisfies

\[
\pi_1\circ h=f,
\qquad
\pi_2\circ h=g.
\]

Thus \texttt{ind2sub} and \texttt{sub2ind} implement the universal property of products.

\subsection*{Coproducts}

The coproduct of two objects is

\[
n_1+n_2.
\]

The canonical injections are

\begin{verbatim}
iota1 = (1:n1)';
iota2 = n1 + (1:n2)';
\end{verbatim}
representing

\[
\iota_1:n_1\to n_1+n_2,
\qquad
\iota_2:n_2\to n_1+n_2.
\]

Given morphisms

\[
f:n_1\to n_B,
\qquad
g:n_2\to n_B,
\]
their copairing is obtained simply by concatenation:

\begin{verbatim}
h = [f;g];
\end{verbatim}
and satisfies

\[
h\circ\iota_1=f,
\qquad
h\circ\iota_2=g.
\]

Thus vector concatenation realizes coproducts.

\subsection*{Pullbacks}

Given morphisms

\[
u:n_A\to n_C,
\qquad
v:n_B\to n_C,
\]
their pullback is obtained by selecting those pairs

\[
(i,j)\in n_A\times n_B
\]
for which

\[
u(i)=v(j).
\]

Using the product construction,

\begin{verbatim}
[pi1,pi2] = ind2sub([nA,nB],(1:nA*nB)');
k = (u(pi1) == v(pi2));
\end{verbatim}

The logical vector \(k\) determines the pullback object as a subobject of the product

\[
n_A\times n_B,
\]
while

\[
\pi_1(k)
\qquad\text{and}\qquad
\pi_2(k)
\]
are the projection morphisms.

A particularly important special case occurs when \(v\) is a monomorphism. Then the pullback computes the inverse image of a subobject and may be implemented directly through

\begin{verbatim}
[k,p] = ismember(u,v);
\end{verbatim}
which explains the central role played by \texttt{ismember} throughout the paper. In that case the canonical projections are obtained as \texttt{pi1=find(k)} and \texttt{pi2=p(k)}.

\subsection*{Equalizers}

Given parallel morphisms

\[
u,v:n_A\to n_B,
\]
the equalizer consists of those elements of \(n_A\) on which the two maps agree.

A direct implementation is

\begin{verbatim}
k = (u == v);
e = find(k);
\end{verbatim}

The vector \(e\) determines the inclusion

\[
e:E\hookrightarrow n_A
\]
of the equalizer object.

A particularly convenient implementation is
\begin{verbatim}
[k,p] = eq(u,v);
\end{verbatim}
with

\begin{verbatim}
k = (u == v);
p = cumsum(k);
\end{verbatim}

Here \(k\) is the characteristic function of the equalizer, while \(p\) records the factorization through the equalizer inclusion. Consequently, every morphism \[ f:n_Z\to n_A \] satisfying \[ u(f)=v(f) \] factors uniquely through the equalizer. More precisely, \[ f=e(\bar f), \] where \[ e=\texttt{find}(k) \] is the equalizer inclusion and \[ \bar f=p(f). \] In Octave notation this becomes 
\begin{verbatim} isequal(f,find(k)(p(f))) \end{verbatim} which returns \texttt{true}. For example,

\begin{verbatim} 
u = [1 2 3 4 5 6]; 
v = [2 2 2 5 5 5];
[k,p] = eq(u,v);
f = [2 2 5 5 2]; isequal(f,find(k)(p(f))) 
\end{verbatim} demonstrates the universal property computationally. The vector \texttt{find(k)} provides the equalizer inclusion, while \(p\) records the unique factorization of \(f\) through that inclusion.

\subsection*{Image Factorizations}

Every morphism

\[
f:n_A\to n_B
\]
admits a canonical factorization into a surjection followed by an injection.

It is computed through

\begin{verbatim}
[Q,r,q] = unique(f);
\end{verbatim}

The vector \(q\) is a surjection

\[
q:n_A\to Q,
\]
while

\[
r:Q\to n_A
\]
is a section satisfying

\[
qr=1_Q.
\]

The original morphism factors as

\[
n_A
\xrightarrow{q}
Q
\xrightarrow{f(r)}
n_B.
\]

Thus

\[
f=(f(r))q
\]
is a (split epimorphism)--monomorphism factorization. In particular, \(Q\) may be identified with the image of \(f\).

\subsection*{Split coequalizers of reflexive pairs}\label{subsec: coeq}

Given a parallel pair

\[
d,c:E\rightrightarrows n,
\]
one first constructs the relation

\begin{verbatim}
R = sparse(c,d,true,n,n);
\end{verbatim}

The smallest equivalence relation containing the identifications imposed by \(d\) and \(c\) is then obtained by computing

\[
(R\cup R^{-1}\cup I)^*.
\]

Using the closure algorithms discussed elsewhere in the paper, this becomes

\begin{verbatim}
T = closure(R,'rst');
\end{verbatim}
or equivalently

\begin{verbatim}
T = transitive(R + R' + I);
\end{verbatim}

The split coequalizer $t$ is constructed as follows:

\begin{verbatim}
[p1,p2] = find(triu(T));
t(p1) = p2;
\end{verbatim}

Since \texttt{find} scans the upper-triangular part in sorted order, the
vector \(t\) assigns to every element the chosen representative of its
equivalence class. Consequently,

\[
t:n\to n
\]
is an idempotent endomorphism satisfying

\[
td=t c=t,
\]
and therefore determines a split coequalizer of the pair

\[
d,c:E\rightrightarrows n.
\]

To recover the corresponding split epimorphism, observe that the fixed
points of \(t\) are precisely the chosen representatives. Thus:

\begin{verbatim}
n=max(d); % d is assumed surjective in Idx

R=logical(sparse(c,d,1,n,n));
T=transitive(R+R'+eye(n));

[p1,p2]=find(triu(T)); % p2 is sorted

t(p1)=p2;  % split coequalizer as an endomorphism of n

r=(t(:)==(1:n)');

q=(cumsum(r))(t);

isequal(t(:),find(r)(q))
\end{verbatim}

The logical vector

\[
r=(t(i)=i)
\]
is the characteristic function of the set of representatives. The map

\[
e=\texttt{find}(r)
\]
is therefore the inclusion

\[
e:Q\hookrightarrow n
\]
of the quotient object \(Q\).

The vector

\[
q=(\texttt{cumsum}(r))(t)
\]
is the corresponding quotient map

\[
q:n\to Q,
\]
and the equality

\begin{verbatim}
isequal(t(:),find(r)(q))
\end{verbatim}
expresses the factorization

\[
t=e\,q.
\]

Since

\[
qe=1_Q,
\]
the pair

\[
q:n\to Q,
\qquad
e:Q\to n
\]
forms a split epimorphism--monomorphism factorization of the split
coequalizer \(t\). Thus the coequalizer construction yields not only the
quotient map but also a canonical choice of representatives for the
equivalence classes determined by the relation generated by \(d\) and
\(c\).

\section{A Boolean Algebra on the Dyadic Interval}\label{Appendix Bollean}

One of the recurring themes of this paper is that infinite mathematical structures can often be represented computationally through finite approximations. An instructive example is provided by the Boolean algebra

\[
\{0,1\}^{\mathbb N},
\]
whose elements may be regarded as infinite binary strings. Since ordinary floating-point numbers already store finite binary expansions, it is natural to use them as coordinates for a computational model of this algebra.

Recall that IEEE 754 double precision represents real numbers using a 52-bit mantissa together with an implicit leading bit. Consequently, the interval

\[
[0.5,1)
\]
may be viewed as a finite approximation of the space of infinite binary strings. Indeed, every number

\[
x\in[0.5,1)
\]
admits a binary expansion

\[
x=0.1b_1b_2b_3\cdots,
\]
and the first 52 binary digits may be recovered exactly.

This correspondence is obtained through the map

\[
x
\longmapsto
x^\sharp
=
\bigl\lfloor (x-0.5)\,2^{53}\bigr\rfloor,
\]
which identifies \(x\) with an integer in the range

\[
0,\ldots,2^{53}-1.
\]

The binary expansion of \(x^\sharp\) then encodes the first 53 bits following the initial digit \(0.1\).

Boolean operations may therefore be transferred from integers to the interval. Given

\[
x,y\in[0.5,1),
\]
define

\[
x\wedge y
=
0.5+
\frac{\operatorname{bitand}(x^\sharp,y^\sharp)}
     {2^{53}},
\]
and

\[
x\vee y
=
0.5+
\frac{\operatorname{bitor}(x^\sharp,y^\sharp)}
     {2^{53}}.
\]

Similarly, the complement operation is given by

\[
\neg x
=
0.5+
\frac{\operatorname{bitcmp}(x^\sharp,53)}
     {2^{53}},
\]
where only the first 53 binary digits are complemented.

Thus the interval acquires the structure of a finite Boolean algebra isomorphic to the algebra of binary strings of length \(53\).

A further refinement is inspired by the familiar ambiguity of dyadic expansions. Every dyadic rational admits two binary representations, one terminating and one eventually repeating:

\[
0.1000\cdots
=
0.0111\cdots.
\]

To retain this information geometrically, one may duplicate each dyadic point and represent its two binary expansions by a pair of conjugate points lying infinitesimally above and below the real axis. The resulting picture resembles a split interval whose non-dyadic points remain on the real line while dyadic rationals bifurcate into conjugate pairs.

This construction is not intended as a rigorous completion of the infinite Boolean algebra \(\{0,1\}^{\mathbb N}\). Rather, it provides a computational model that simultaneously preserves its binary structure, its geometric intuition, and its compatibility with standard floating-point arithmetic. It also illustrates a recurring principle of computational category theory: infinite mathematical objects are often best understood through finite coordinate systems that retain the essential algebraic structure while remaining directly computable.

\stop
\begin{lstlisting}[language=Octave,
caption={Visualization of the split dyadic interval}]
% fast_dyadic_boolean_plot.m
%
% Visualization of the dyadic interval viewed as a
% finite approximation of the Boolean algebra
% {0,1}^N.
%
% Each dyadic rational is represented by two
% conjugate points corresponding to its terminating
% and repeating binary expansions.

clear;
clc;
clf;

% Maximum dyadic degree 2^{-k}
max_degree = 8;

% Vertical separation of conjugate pairs
epsilon_base = 0.1;

real_parts = [];
imag_parts = [];
colors = [];

% Generate dyadic rationals by degree

for k = 1:max_degree

    denom = 2^k;

    % Only odd numerators are used in order
    % to generate each dyadic exactly once.

    for m = 1:2:denom

        val = m / denom;

        % Dyadics born at deeper levels are drawn
        % closer to the real axis.

        h = epsilon_base / k;

        % Terminating expansion

        real_parts(end+1) = val;
        imag_parts(end+1) = h;
        colors(end+1,:) = [0.2, 0.6, 0.9];

        % Repeating expansion

        real_parts(end+1) = val;
        imag_parts(end+1) = -h;
        colors(end+1,:) = [0.9, 0.3, 0.3];

    end

end

figure(1);
hold on;

% Real axis

plot([0,1],[0,0], ...
     'k--', ...
     'LineWidth',1.5);

% Conjugate dyadic representatives

scatter(real_parts, ...
        imag_parts, ...
        25, ...
        colors, ...
        'filled');

title(sprintf( ...
    'Split Dyadic Interval (Max Degree: 2^{-%d})', ...
    max_degree), ...
    'FontSize',12);

xlabel('Dyadic Value', ...
       'FontSize',11);

ylabel('Conjugate Height', ...
       'FontSize',11);

xlim([-0.05,1.05]);
ylim([-epsilon_base*1.5,epsilon_base*1.5]);

grid on;
box on;

hold off;
\end{lstlisting}

Dyadic rationals are generated exactly once by considering only odd numerators. Each dyadic point is subsequently replaced by two representatives corresponding to its terminating and repeating binary expansions, drawn respectively above and below the real axis. The resulting picture provides a geometric realization of the interval as a Boolean-algebraic space of infinite binary strings, with dyadic points resolved into conjugate pairs.

\end{document}